\documentclass[11pt]{amsart}
\usepackage{amsmath,amssymb,amsthm,mathtools}
\usepackage[margin=1.1in]{geometry}
\usepackage{enumitem}
\usepackage{iftex}
\ifPDFTeX\else\usepackage[no-math]{fontspec}\fi
\usepackage{hyperref}
\usepackage{cleveref}
\newtheorem{theorem}{Theorem}[section]

\newtheorem{corollary}[theorem]{Corollary}
\theoremstyle{definition}
\newtheorem{definition}[theorem]{Definition}
\theoremstyle{remark}
\newtheorem{remark}[theorem]{Remark}
\newcommand{\eps}{\varepsilon}
\newcommand{\Pplus}{P^{+}}

\newcommand{\repo}{\url{https://github.com/daizisheng/erdos-824}}

\title[Coprime pairs with equal sum of divisors]{Coprime pairs with equal sum of divisors,\\ assuming smooth shifted primes}
\author{Shisheng Li}
\author{Zhong Zhuang}
\date{Draft, 11 October 2026.}

\begin{document}
\sloppy
\begin{abstract}
Let $h(x)$ be the number of pairs $1\le a<b<x$ with $\gcd(a,b)=1$ and $\sigma(a)=\sigma(b)$. Erd\H{o}s Problem \#824 asks whether $h(x)>x^{2-o(1)}$. By an argument of Pollack and Pomerance, this follows if for every fixed $\delta>0$ there are $y^{1-o(1)}$ primes $p\le5y$ with $p+1$ free of prime factors exceeding $y^\delta$. We write this implication out with explicit parameters and formalize it in Lean~4. We then check that the argument of the recent OpenAI preprint \emph{Weighted dilation graphs, smooth shifted primes and totient fibers}, which claims the corresponding statement for $p-1$, applies to $p+1$ after a change of sign in two of its sections. Consequently $h(x)=x^{2-o(1)}$, \emph{conditionally on the correctness of that preprint}, which has not been refereed and whose main new estimates we have not independently verified.
\end{abstract}

\maketitle

\section{Introduction}

Let $\sigma(n)$ be the sum of the divisors of $n$ and
\[
 h(x)=\#\{(a,b):\ 1\le a<b<x,\ \gcd(a,b)=1,\ \sigma(a)=\sigma(b)\}.
\]
Trivially $h(x)<x^2$. Erd\H{o}s~\cite{Er59} asked whether $h(x)/x\to\infty$; Problem \#824 of~\cite{EP} asks whether $h(x)>x^{2-o(1)}$. Pollack and Pomerance~\cite[\S6]{PP16} answered Erd\H{o}s's question, proving $h(x)>x^{1.4}$ for large $x$. Following an idea of Erd\H{o}s~\cite{Er35}, their argument shows: if $\theta>1$ is such that there are $y^\theta/(\log y)^{O(1)}$ primes $p\le y^\theta$ with $\Pplus(p+1)\le y$, then there are at least $x^{(2-\eps)(1-1/\theta)+o(1)}$ coprime pairs, for every fixed $\eps>0$. Here $\Pplus(m)$ is the largest prime factor of $m$. With the Baker--Harman value $\theta\approx3.377$~\cite{BH98} this gives the exponent $1.4$; an admissible $\theta$ for every $\theta>1$ would give $x^{2-o(1)}$, as also remarked by Tao on the problem's forum~\cite{Tao}. A claim of $h(x)\ge x^{13/8-o(1)}$ by a different route was submitted to the forum in August 2026~\cite{Cam}.

For a fixed smoothness exponent, the best unconditional result is due to Lichtman~\cite{Li22}: for fixed $a\ne0$ and $\delta>15/(32\sqrt e)=0.2843\ldots$ there are $\gg x/(\log x)^C$ primes $x<p\le2x$ with $\Pplus(p-a)\le x^\delta$. The corresponding statement for every fixed $\delta>0$ is a conjecture of Erd\H{o}s. For $a=1$ it is claimed in the OpenAI preprint~\cite{WD} (``WD'' below):

\begin{theorem}[{\cite[Theorem 1.2]{WD}}]\label{thm:wd}
For every fixed $0<\delta<1/4$, as $x\to\infty$,
$\#\{p\text{ prime}:\ 2x<p\le5x,\ \Pplus(p-1)\le x^\delta\}\ge x^{1-o(1)}$.
\end{theorem}

The application to $\sigma$ needs $p+1$ instead of $p-1$.

\begin{definition}\label{def:aplus}
Hypothesis $\mathrm{A}^+$: for every $0<\delta<1/4$ and every $\eta>0$ there is $y_0\ge1$ such that for every real $y\ge y_0$,
\[
 \#\{p\text{ prime}:\ p\le 5y,\ \Pplus(p+1)\le y^\delta\}\ \ge\ y^{1-\eta}.
\]
\end{definition}

\begin{theorem}\label{thm:main}
If $\mathrm{A}^+$ holds, then for every $\eps>0$ we have $h(x)>x^{2-\eps}$ for all sufficiently large $x$. The pairs may moreover be taken with $a,b$ squarefree.
\end{theorem}

\Cref{thm:main} is the argument of~\cite[\S6]{PP16}, with $\theta$ chosen large but fixed in terms of $\eps$. We give it in full in \S\ref{sec:reduction}, with explicit parameters and under the weaker hypothesis of $y^{1-o(1)}$ (rather than $y/(\log y)^{O(1)}$) smooth shifted primes, because this is the version we have formalized in Lean~4 (\S\ref{sec:lean}).

\begin{theorem}\label{thm:aplus}
Assume that the argument of~\cite{WD} proving \Cref{thm:wd} is correct. Then the same argument, modified as described in \S\ref{sec:sign}, shows that for every fixed $0<\delta<1/4$, $\#\{p\text{ prime}:\ 2x-2<p\le5x,\ \Pplus(p+1)\le x^\delta\}\ge x^{1-o(1)}$; in particular $\mathrm{A}^+$ holds.
\end{theorem}

\begin{corollary}
If the argument of~\cite{WD} is correct, then $h(x)=x^{2-o(1)}$.
\end{corollary}

The status of the corollary is exactly that of~\cite{WD}, an unrefereed preprint. Its Sections~3--5 introduce a new transference and shifted-correlation method; we have not verified them independently (\S\ref{sec:status}). Our contribution is the observation that the sign change goes through, a check of the sections that it affects, and the formalization of \Cref{thm:main}.

\section{Smooth successors from the argument of WD}\label{sec:sign}

\Cref{thm:aplus} is proved by the argument of~\cite{WD} for \Cref{thm:wd}, which detects primes of the form $2u+1$ with $u$ supported on $x^\delta$-smooth integers; we apply it to primes of the form $2u-1$ instead, so that $p+1=2u$ is $x^\delta$-smooth. The results of Sections~2--5 of~\cite{WD} are used as stated. In Sections~6--7 of~\cite{WD} (the Type II estimate and the prime extraction) one replaces $2u+1$ by $2u-1$ throughout and adjusts the congruences accordingly; the local densities and residue counts of the sieve are the same for both signs.

The one step that is not purely notational is the determinant parametrization in the Type II estimate. Two factorizations $mn=2eh-1$, $m'n=2eh'-1$ with $m\ne m'$ give $m'-m=-2ek$ and $h'-h=-kn$ for some integer $k\ne0$, and then $mh'-m'h=-k(mn-2eh)=k$. So the shifted endpoints are $z=m'h$ and $z+k=mh'$, as in~\cite{WD}, and the shift-cancellation theorem~\cite[Theorem~5.1]{WD}, which holds for all $0<|k|\le(\log x)^C$, applies unchanged; the original argument already uses both signs of $k$. The converse reconstruction of the factorizations and the divisor-moment bounds, now with $\tau(2u-1)$ and $\tau(v+1)$, go through in the same way.

\section{Proof of \Cref{thm:main}}\label{sec:reduction}

It suffices to treat $0<\eps<1$: for $\eps\ge1$ take any $\eps'\in(0,1)$; then $h(x)>x^{2-\eps'}\ge x^{2-\eps}$ for $x\ge1$. Write $L=\log x$. Fix $\eps\in(0,1)$, a real $B\ge8/\eps$, $\eta=\eps/16$ and $\delta=1/(2B)\in(0,1/4)$; these are fixed before $x\to\infty$. For large $x$ put $y=L^B$,
\[
 \mathcal P=\{p\text{ prime}:\ p\le5y,\ \Pplus(p+1)\le y^\delta\},\qquad M=|\mathcal P|.
\]
By $\mathrm{A}^+$, $M\ge y^{1-\eta}=L^{B(1-\eta)}$ for all large $x$. Let $k=\lfloor L/\log(5y+1)\rfloor$. Then $(5y+1)^k\le x$, and since $\log(5y+1)=B\log L+O(1)$,
\begin{equation}\label{eq:k}
 k=\frac{L}{B\log L}(1+o(1)),\qquad k\le L .
\end{equation}
For large $x$, $M\ge L^2\ge2k$ because $B(1-\eta)>2$.

\subsection{The family} Let $\mathcal A$ be the set of products of $k$ distinct primes of $\mathcal P$. Each $a\in\mathcal A$ is squarefree and $a<(5y+1)^k\le x$, and distinct $k$-subsets give distinct products. So $N:=|\mathcal A|=\binom Mk\ge(M/k)^k$ and, by~\eqref{eq:k} and $\log k\le\log L$,
\begin{equation}\label{eq:N}
 \log N\ge k\log L\,(B(1-\eta)-1)=L\Bigl(1-\eta-\frac1B\Bigr)(1+o(1)).
\end{equation}
For $a=p_1\cdots p_k$ we have $\sigma(a)=\prod_i(p_i+1)\le x$, and every prime factor of $\sigma(a)$ is at most $y^\delta=L^{1/2}$. The number $U$ of $L^{1/2}$-smooth integers in $[1,x]$ satisfies
\begin{equation}\label{eq:U}
 U\le\Bigl(1+\frac L{\log2}\Bigr)^{L^{1/2}},\qquad \log U=o(L),
\end{equation}
since such an integer is $\prod_{q\le L^{1/2}}q^{e_q}$ with $0\le e_q\le L/\log2$.

\subsection{Many equal values} With $r(s)=\#\{a\in\mathcal A:\sigma(a)=s\}$, Cauchy's inequality gives
\[
 \#\{(a,b)\in\mathcal A^2:\ \sigma(a)=\sigma(b)\}=\sum_sr(s)^2\ge N^2/U .
\]

\subsection{Pairs sharing many primes} Let $t=\lceil\eta k\rceil$. The ordered pairs in $\mathcal A^2$ sharing at least $t$ primes number at most $N\binom kt\binom M{k-t}$ (choose $a$, $t$ of its primes, and the other $k-t$ primes of $b$), and
\[
 \binom M{k-t}\Big/\binom Mk=\prod_{i=0}^{t-1}\frac{k-i}{M-k+t-i}\le\Bigl(\frac{2k}M\Bigr)^t,\qquad\binom kt\le2^k .
\]
By~\eqref{eq:k}, \eqref{eq:U} and $M\ge L^{B(1-\eta)}$,
\begin{align*}
 \log\bigl(2^k(2k/M)^t\cdot2U\bigr)&\le(k+1)\log2+\log U-\eta k\bigl(B(1-\eta)\log L-\log(2L)\bigr)\\
 &=-\eta L\Bigl(1-\eta-\frac1B\Bigr)(1+o(1))+o(L),
\end{align*}
which is negative for large $x$. So at least $N^2/(2U)$ ordered pairs $(a,b)$ have $\sigma(a)=\sigma(b)$ and share fewer than $t$ primes; in particular $a\ne b$.

\subsection{Removing the common factor} For such a pair let $g=\gcd(a,b)$, $u=a/g$, $v=b/g$. As $a,b$ are squarefree, $g,u,v$ are pairwise coprime, so $\sigma(g)\sigma(u)=\sigma(g)\sigma(v)$ and $\sigma(u)=\sigma(v)$, with $\gcd(u,v)=1$, $u\ne v$, and $u,v<x$ squarefree. The integer $g$ is a product of at most $t-1<\eta k$ primes below $5y+1$, so $g\le(5y+1)^{\eta k}\le x^\eta$. Given the ordered pair $(u,v)$, the pair $(a,b)=(gu,gv)$ is determined by $g$; each unordered pair $\{u,v\}$ therefore arises from at most $2x^\eta$ ordered pairs $(a,b)$. Hence
\[
 h(x)\ge\frac{N^2}{4Ux^\eta},\qquad
 \log h(x)\ge\Bigl(2-3\eta-\frac2B\Bigr)L-o(L)>(2-\eps)L
\]
for large $x$, because $3\eta+2/B\le3\eps/16+\eps/4<\eps/2$. \qed

\begin{remark}
Since $h(x)$ only depends on the integers below $x$, $h(x)=h(\lceil x\rceil)$. Hence the theorem for integer $x$, which is the form formalized in Lean, implies it for real $x$: for $0<\eps<1$, $h(x)=h(\lceil x\rceil)>\lceil x\rceil^{2-\eps}\ge x^{2-\eps}$. Pollack and Pomerance count pairs with $b\le x$; the difference is immaterial.
\end{remark}

\section{Formalization}\label{sec:lean}

The lower bound in \Cref{thm:main} is formalized in Lean~4 (v4.35.0-rc2) with Mathlib (commit \texttt{0653561}) in the directory \texttt{lean/} of~\cite{repo} (main theorem in \texttt{lean/Erdos824/Main.lean}). The definitions are
\begin{verbatim}
noncomputable def h (x : ℕ) : ℕ :=
  ((Finset.Ico 1 x ×ˢ Finset.Ico 1 x).filter
    (fun p => p.1 < p.2 ∧ Nat.Coprime p.1 p.2 ∧ sigma 1 p.1 = sigma 1 p.2)).card
def SmoothSuccessorHypWeak : Prop :=
  ∀ δ : ℝ, 0 < δ → δ < 1/4 → ∀ η : ℝ, 0 < η → ∃ y0 : ℝ, ∀ y : ℝ, y0 ≤ y →
    (y ^ (1 - η) : ℝ) ≤ (((Finset.Iic ⌊5*y⌋₊).filter
      (fun p => p.Prime ∧ ∀ q ∈ (p+1).primeFactors, ((q:ℕ):ℝ) ≤ y ^ δ)).card : ℝ)
\end{verbatim}
where \texttt{sigma 1} is Mathlib's $\sigma_1=\sigma$; \texttt{SmoothSuccessorHypWeak} is \Cref{def:aplus}. The theorem
\begin{verbatim}
theorem erdos824_of_smoothSuccessorWeak (H : SmoothSuccessorHypWeak) :
    ∀ ε : ℝ, 0 < ε → ∃ X : ℝ, ∀ x : ℕ, X ≤ x → (x:ℝ) ^ (2 - ε) < (h x : ℝ)
\end{verbatim}
formalizes the lower bound of \Cref{thm:main} (not the squarefree addendum, which is proved only in \S\ref{sec:reduction}). It has no \texttt{sorry}, and \texttt{\#print axioms} reports only \texttt{propext}, \texttt{Classical.choice} and \texttt{Quot.sound}. The Lean proof follows \S\ref{sec:reduction} with $B=16/\eps$ and $\delta=\eps/32$.

\section{Status}\label{sec:status}

The lower-bound implication in \Cref{thm:main} is unconditional and machine-checked; the squarefree addendum is not part of the formalization. \Cref{thm:aplus} rests on~\cite{WD}. We rechecked Sections~6--7 of WD for the change of sign, but we have not independently verified Sections~3--5 (transference for weighted dilation graphs, construction of the ideal operator, lifted shift cancellation; about 2800 lines of source). Writing Lean statements for WD raised a hypothesis clarification: the definition of prime groups writes $\ell\in\mathbb N$, whereas the proof in Section~4 uses $\ell\ge1$. The application in Section~7 takes $\ell=1$, so this point does not obstruct the present application.

\section*{Acknowledgements}
The reduction in \S\ref{sec:reduction} is due to Pollack and Pomerance~\cite{PP16}, after Erd\H{o}s; its relevance to Problem \#824 was also pointed out by T.~Tao. AI tools were used extensively in this work. The project owner reports that ChatGPT (OpenAI) produced the Lean formalization and that it has not undergone human mathematical review. Codex (OpenAI) adapted the formalization to the current Lean toolchain and Palomar submission format. Mechanical verification of the conditional implication does not constitute human review or verify Hypothesis~$\mathrm{A}^+$ or the argument of~\cite{WD}. The preprint~\cite{WD} is itself AI-generated.

\end{document}